\documentclass[12pt]{amsart}
\usepackage{amscd,times,fullpage}
\usepackage[utf8x]{inputenc}

\usepackage{amsthm,amsmath}
\usepackage{amssymb}
\usepackage{graphicx}
\usepackage{amsthm}
\usepackage{enumerate}
\usepackage{dsfont}
\usepackage{tikz-cd}
\usetikzlibrary{cd}
\usetikzlibrary{positioning}
\usetikzlibrary{matrix}
\usepackage{hyperref}

\hypersetup {
	pdfpagemode = {UseNone},
	pdftitle = {Embeddings of Homogeneous Sasakian Manifolds},
	pdfauthor = {Andrea Loi, Giovanni Placini, Michela Zedda},
	pdflang = {en-UK}
} 
\usepackage{thmtools}
\usepackage{thm-restate}

\usepackage{cleveref}

\newcommand{\R}{\mathbb{R}}

\newcommand{\Z}{\mathbb{Z}}

\newcommand{\CP}{\mathbb{C}\mathrm{P}}
\newcommand{\N}{\mathbb{N}}

\newcommand{\CH}{\mathbb{C}\mathrm{H}}

\newcommand{\C}{\mathbb{C}}            
\newcommand{\de}{\partial}

\newcommand{\K}{K\"{a}hler }

\newcommand{\OO}{\mathcal{O}}

\newcommand{\ov}[1]{\overline{#1}}

\newcommand{\deb}{\ov\partial}

\newcommand{\di}{{\operatorname{d}}}

\newcommand{\Id}{\operatorname{Id}}

\newcommand{\lra}{\longrightarrow}

\newcommand{\D}{\mathcal{D}}

\newcommand{\TM}{\text{T} M}

\newtheorem{theor}{Theorem}[section]
\newtheorem{prop}[theor]{Proposition}
\newtheorem{defin}[theor]{Definition}
\newtheorem{lem}[theor]{Lemma}
\newtheorem{cor}[theor]{Corollary}

\newtheorem{ex}[theor]{Example}
\newtheorem{remark}[theor]{Remark}

\usepackage{nicefrac}

\begin{document}

	\title[Immersions into Sasakian space forms]{Immersions into Sasakian space forms}
	
	\author{A.~Loi}
	\address{Dipartimento di Matematica e Informatica, Universit\'a degli studi di Cagliari, Via Ospedale 72, 09124 Cagliari, Italy}
	\email{loi@unica.it}
	\author{G.~Placini}
	\address{Dipartimento di Matematica e Informatica, Universit\'a degli studi di Cagliari, Via Ospedale 72, 09124 Cagliari, Italy}
	\email{giovanni.placini@unica.it}
	\author{M.~Zedda}
	\address{Dipartimento di Scienze Matematiche, Fisiche e Informatiche, Universit\'a di Parma, Parma, Italy}
	\email{michela.zedda@unipr.it}

	\date{\today ; {\copyright  A.~Loi, G.~Placini and M.~Zedda 2023}}
	
	\subjclass[2010]{53C25; 53C30; 53C24; 53C42; 53A55}
	\keywords{Sasakian space form; Sasakian immersion, Calabi Rigidity, Homogeneous Sasakian manifold}
	\thanks{The authors are supported by INdAM and  GNSAGA - Gruppo Nazionale per le Strutture Algebriche, Geometriche e le loro Applicazioni. The first two authors are supported by GOACT - Funded by Fondazione di Sardegna and funded by PNRR e.INS Ecosystem of Innovation for Next Generation Sardinia (CUP F53C22000430001, codice MUR ECS00000038).}

	\begin{abstract}
		We study immersions of Sasakian manifolds into finite and infinite dimensional Sasakian space forms. 
		After proving Calabi's rigidity results in the Sasakian setting, we characterise all homogeneous Sasakian manifolds which admit a (local) Sasakian immersion into a nonelliptic Sasakian space form.
		Moreover, we give a characterisation of homogeneous Sasakian manifolds which can be embedded into the standard sphere both in the compact and noncompact case.
		
	\end{abstract}
	
	\maketitle

	\section{Introduction and statements of the main results}\label{sectionint}
	
	A Sasakian manifold $M=(M,\eta,R,g,\Phi)$ is a contact manifold equipped with a compatible Riemannian metric and a transverse complex structure defining a \K structure transverse to the Reeb foliation. 
	Other than the transverse one, there is another natural \K structure associated to a Sasakian manifold: the one on its Riemannian cone.
	Due to their deep connection, many problems in \K geometry can be studied in the Sasakian setting with analogous results.
	This paper provides such an instance.
	Namely, we are interested in Sasakian manifolds admitting a Sasakian immersion or embedding in a Sasakian space form.
	Although this problem has a long tradition in the \K case initiated in the 50's by Calabi~\cite{Calabi53Isometric} (see~\cite{Loi18Book} for a recent review), the Sasakian counterpart has recently received relatively little attention.
	This problem has been studied in the late 60's and 70's under different names, see e.g. Harada~\cite{Harada72SpaceFormInSpaceForm,Harada72SubmanifoldsI,Harada72SubmanifoldsII,Harada75SomeSasakianSubmanifolds}, Kenmotsu~\cite{Kenmotsu69InvariantSubmanifolds}, Kon~\cite{Kon73Kodai,Kon76InvariantSubmanifolds}, Okumura~\cite{Okumura68Immersion}, and only more recently by Bande, Cappelletti--Montano and the first two named authors~\cite{Bande20EtaEinsteinNonCompact,Cappelletti19EinsteinSpheres,Placini21SasakiSolitons}.
	
	We are interested in Sasakian immersions into simply connected Sasakian space forms and we initially focus on Sasakian analogues of some classical results of Calabi. 
	As in the \K case, given $N\leq\infty$, there are only three complete simply connected Sasakian space forms $M(N,c)$ (up to transverse deformation) of dimension $2N+1$ distinguished by the $\phi$-sectional curvature $c$. Namely, the hyperbolic space form with $c<-3$, the Heisenberg space with $c=-3$ and the standard sphere with $c>-3$. 
	Moreover, Sasakian space forms are bundles over complex space forms, trivial ones in the first two cases. 
	For this reason, when the Sasakian manifolds considered are regular, some statements follow directly from the \K case using the fact that a Sasakian immersion induces a \K immersion of the space of Reeb leaves.
	As an example, the classification of immersions of Sasakian space forms into Sasakian space forms was already pointed out by Harada~\cite{Harada75SomeSasakianSubmanifolds} without claims of originality.
	We exploit the fibration of Sasakian space forms repeatedly in this paper. 
	For instance, the proof of our first result, the Sasakian analogue of Calabi's rigidity, relies on this fact.
	\begin{restatable}[Sasakian rigidity]{theor}{ThmSasakianRigidity}
		\label{ThmSasakianRigidity}
		Let $\varphi_1,\varphi_2:M\lra M(n,c)$ be Sasakian immersions of a Sasakian manifold $M$ into a Sasakian space form $M(N,c)$ with $N\leq\infty$. Then there exists a Sasakian transformation $T$ of $M(N,c)$ such that $\varphi_1=T\circ\varphi_2$.
	\end{restatable}
	Notice that this is a peculiarity of space forms. 
	In fact, one can find many examples of Sasakian manifolds which admit immersions in a Sasakian manifold $\hat{M}$ not related by an isometry of $\hat{M}$.
	For instance, one can take a $\hat{M}$ to be a regular Sasakian manifold constructed over the examples provided by Green~\cite{Green78NonRigidity}.
	Another classical result, whose proof follows \textit{mutatis mutandis} Calabi's proof, is an extension result for local immersions.
	\begin{restatable}[Global extension]{theor}{ThmSasakianGluing}
		\label{ThmSasakianGluing}
		Let $M$ be a simply connected Sasakian manifold. If for all $p\in M$ there exists a local Sasakian immersion $\varphi_p:U_p\lra M(N,c)$ of  a neighbourhood $U_p$ of $p$  into a Sasakian space form $M(N,c)$ with $N\leq\infty$, then $M$ admits a global immersion into $M(N,c)$.
		Furthermore, the immersion is unique up to rigid transformations of $M(N,c)$.
	\end{restatable}
	This is clearly false for manifolds with nontrivial fundamental group. For example, the compact Heisenberg manifold, being the quotient of $\R^{2n+1}=M(n,-3)$ by the integral Heisenberg group, admits a local (but not global) immersion into $M(n,-3)$ itself.
	The hypothesis of existence of a local immersion at each point can be soften to the existence of a local immersion at one point in some peculiar cases, e.g. when the Sasakian manifold is locally homogeneous or the metric is real analytic.
	
	Given a Sasakian manifold $M$, deciding whether $M$ admits a Sasakian embedding or immersion into a Sasakian space forms is a highly nontrivial problem.
	In light of this we begin by investigating the special case of (locally) homogeneous Sasakian manifolds. 
	In this case, similarly to the \K setting, there is a dichotomy. 
	Namely, immersions into Sasakian space forms of $\phi$-sectional curvature $c\leq-3$ are much more restricted than immersions into spheres.
	In the Sasakian setting this is already noticeable in \cite{Bande20EtaEinsteinNonCompact} where the authors study local immersions of $\eta$-Einstein manifolds into Sasakian space forms of $\phi$-sectional curvature $c\leq-3$.
	In fact, we obtain a similar result for local immersions of locally homogeneous Sasakian manifolds which is the Sasakian analogue of a theorem of Di~Scala, Ishi and Loi~\cite{Discala12Immersions}. 
	Namely, writing $F(n,b)$ for the complex space form of dimension $2n$ and holomorphic sectional curvature $4b$, we prove the following
	\begin{restatable}{theor}{mainHom}
		\label{mainHom}
		Let $M$ be a complete locally homogeneous Sasakian manifold of dimension $2n+1$. Suppose there exists an open subset $U\subset M$ and a Sasakian immersion $\varphi:U\lra M(N,c)$ with $N\leq\infty$.  
		\begin{enumerate}
			\item If $c<-3$, $M$ is Sasaki equivalent to $M(n,c)/\Gamma$ for a discrete group of Sasakian transformations of $M(n,c)$. Moreover, if $U=M$, then $\Gamma=\{1\}$. \label{point1mainHom}
			\item if $c=-3$, then $M$ is Sasaki equivalent to the quotient of 
			$$\R\times F(k,0)\times F(n_1,b_1)\times\cdots\times F(n_r,b_r)$$
			with $k+n_1+\dots+n_r=n$, and  $b_i<0$ for $i=1,\dots, r$,
			by a discrete group of Sasakian transformations $\Gamma$.
			Moreover, if $U=M$, the group $\Gamma$ is trivial.\label{point2mainHom}
		\end{enumerate}
	\end{restatable}
	Let us briefly comment on this result. 
	First of all, notice that any of the above Sasakian manifolds with nontrivial fundamental group shows the necessity of the assumption on $\pi_1(M)$ in Theorem~\ref{ThmSasakianGluing}.
	One can describe the fundamental groups $\Gamma=\pi_1(M)$ explicitly. 
	More specifically, $\Gamma$ is the extension of a product of discrete subgroups of \K isometries of the factors $F(n_i,b_i)$, cf. Remark~\ref{RmkFundGroup}.
	When the immersion considered in Theorem~\ref{mainHom} is global we can describe the immersion explicitly, see Remark~\ref{RmkGlobalImmersion}.
	In particular, if we assume $N<\infty$, the immersion is just the linear inclusion of $M(n,c)$ in $M(N,c)$.
	Observe that the same thesis is obtained in \cite{Bande20EtaEinsteinNonCompact} under the hypotheses that $M$ is $\eta$-Einstein and $N<\infty$. 
	Regarding infinite dimensional Sasakian space forms, note that $M(n,c)$ is $\eta$-Einstein and admits an immersion in $M(\infty,-3)$ for $c<-3$ by Theorem~\ref{mainHom}. 
	An interesting open problem, closely related to a conjecture in the \K setting~\cite[Conjecture~4.1.4]{Loi18Book}, is to characterise $\eta$-Einstein manifolds which admit an immersion into $M(\infty,-3)$.
	
	In analogy with the \K realm, Sasakian immersions in space forms $M(N,c)$ with c$>-3$ are much more flexible. 
	This is ultimately a consequence of Kodaira Embedding Theorem and its Sasakian version explained in detail in Section~\ref{SecGeneralSasakianEmbedding} below.
	Namely, the pullback of the Hopf bundle $M(N,1)\lra F(N,1)$ via the Kodaira embedding $X\lra F(N,1)$ of any projective manifold $X$ provides a regular Sasakian manifold $M$ together with a Sasakian embedding $M\lra M(N,1)$.
	On the other hand, given a Sasakian manifold $M$, it is interesting to know whether the Sasakian structure on $M$ is of this type, that is, if $M$ admits an embedding into a standard sphere.
	By the structure Theorem~\ref{TheorStructure} a regular Sasakian manifold $M$ is, up to transverse homothety, the unitary bundle of a hermitian line bundle $(L^{-1},h^{-1})$ over a \K manifold $(X,\omega)$ with $c_1(L)=[\omega]$.
	Therefore, the problem at hand is strictly related to finding an embedding of $(X,L)$ into $(\CP^N, \OO(1))$ isometric for the hermitian metrics on the line bundles.

	As in the nonelliptic case, we focus on homogeneous manifolds for the initial approach to the problem.
	Notice that if we restrict to finite dimensional case, a Sasakian manifold admitting an immersion in $M(N,1)$ is forced to be compact, see Remark~\ref{infinitenoncompact}.
	Our next result shows that for homogeneous manifolds this condition is also sufficient and the immersion is in fact an embedding. Moreover, we are able to describe the embedding explicitly.
	\begin{restatable}{theor}{TheoMainCpt}\label{TheoMainCpt}
		Up to transverse homothety, a compact homogeneous Sasakian manifold $M$ admits a Sasakian embedding $\varphi$ into a finite dimensional standard sphere $S^{2N+1}$ which is unique up to unitary transformations of $\C^{N+1}$.			
		
		Moreover, identifying $M$ with the unitary bundle associated to the line bundle $(L^{-1},h^{-1})$ over a \K manifold $X$, the embedding is given by
		\begin{align*}
			\varphi: M\hspace{1mm}&\longrightarrow \hspace{8mm}S^{2N+1}\\
			v\ \ &\mapsto \dfrac{\left(v(s_0(x)), \dots , v(s_{N}(x))\right)}{\sum_{j=0}^N h\left(s_j(x),s_j(x)\right)}
		\end{align*}
		where  $v$ lies in the fibre of $L^{-1}$ over $x\in X$ and $(s_0,\ldots,s_N)$ is a suitable basis of $H^0(L)$.
	\end{restatable}
	It was proven in~\cite{Takeuchi78Homogeneous} that a \K manifold admitting an immersion in $\CP^N$ is necessarily simply connected. 
	Although this is not the case in the Sasakian setting, the absence of an assumption on the fundamental group is not unexpected because a compact homogeneous Sasakian manifold necessarily has cyclic fundamental group.
	One should compare this result with~\cite{Cappelletti19EinsteinSpheres} where Cappelletti--Montano and the first named author characterised $\eta$-Einstein manifolds admitting a codimension $2$ immersion in $M(N,1)$. 
	
	If we allow the standard sphere to be infinite dimensional, there are plenty of Sasakian manifolds which can be immersed into $M(N,1)$.
	Even restricting to the homogeneous case, one can see that compactness is not a necessary condition.
	In the \K case, up to homotheties, a homogeneous manifold admits an immersion into $\CP^\infty$ if and only if it is simply connected, see~\cite[Theorem~3]{Discala12Immersions} and \cite[Theorem~1.1]{Loi15HomogeneousImmersions}.
	We prove that the characterisation of homogeneous Sasakian manifolds which can be immersed in the infinite dimensional sphere is analogous.
	\begin{restatable}{theor}{TheoMainTop}\label{TheoMainTop}
		Let $M$ be a homogeneous Sasakian manifold. Then a suitable $\D$-homothetic transformation of $M$ admits a Sasakian immersion $\varphi:M\lra S^\infty$ if and only if $\pi_1(M)$ is cyclic.
	\end{restatable}
	While in the compact case a transverse homothety is necessary to make the transverse \K form integral, here this condition is essentially different.
	In fact, there are contractible homogeneous Sasakian manifolds which cannot be immersed into $S^\infty$ even though every transverse \K class is trivial and thus integral, see Remark~\ref{necessaryDhomothetic}.
	Observe that the immersion in Theorem~\ref{TheoMainTop} is full, i.e. its image is not contained in a proper totally geodesic submanifold, if and only if the manifold $M$ is noncompact.
	This is no longer true if we drop the homogeneity assumption, even if the metric is complete and $N<\infty$. See Example~\ref{ExCompleteFiniteInjective} for an injective Sasakian immersion of a noncompact complete Sasakian manifold into a finite dimensional standard sphere.
	Putting these together we get
	\begin{cor}
		Up to transverse homothety, a homogeneous Sasakian manifold $M$ admits a Sasakian embedding $\varphi$ into a finite dimensional standard sphere $S^{2N+1}$ if and only if it is compact.		
	\end{cor}

	\noindent \textbf{Structure of the paper.} The paper is organised as follows.  In  Section~\ref{sectionbackground} we review the basics of Sasakian geometry with particular focus on homogeneous Sasakian manifolds, Sasakian immersions and regular Sasakian structures.
	The main results of this paper are divided into three sections. 
	Namely, in Section~\ref{SecRigidity} we discuss immersions of Sasakian manifolds in Sasakian space forms and prove Sasakian rigidity (Theorem~\ref{ThmSasakianRigidity}) and the global extension theorem (Theorem~\ref{ThmSasakianGluing}).
	Section~\ref{SecImmersionNegative} is dedicated to the proof of Theorem~\ref{mainHom}.
	The remainder of the paper deals with immersions into finite and infinite dimensional spheres.
	In particular, in Section~\ref{SecGeneralSasakianEmbedding} we discuss CR embeddings of regular Sasakian manifolds into standard spheres and discuss the induced Sasakian structure.
	This discussion is used in Section~\ref{SecImmersionPositive} to prove Theorem~\ref{TheoMainCpt} and Theorem~\ref{TheoMainTop}.

	\section{Sasakian manifolds}\label{sectionbackground}

	Sasakian geometry can be understood in terms of contact metric geometry and via the associated \K cone, cf. the monograph of Boyer and Galicki \cite{Boyer08Book}. We will present both formulations for the reader convenience, but we will focus mostly on the regular case for it is central in this paper.
	In the following all manifolds and orbifolds are assumed to be connected.

	A \textit{K-contact structure} $(\eta,\Phi,R,g)$ on a manifold $M$ consists of a contact form $\eta$ and an endomorphism 
	$\Phi$ of the tangent bundle $\TM$ satisfying the following properties:
	\begin{enumerate}
		\item[$\bullet$] $\Phi^2=-\Id+R\otimes\eta$ where $R$ is the Reeb vector field of $\eta$,
		\item[$\bullet$] $\Phi_{\vert\D}$ is an almost complex structure compatible with the symplectic form $\di\eta$ on $\D=\ker\eta$,		
		\item[$\bullet$] the Reeb vector field $R$ is Killing with respect to the metric $g(\cdot,\cdot)=\dfrac{1}{2}\di\eta(\cdot,\Phi\cdot)+\eta(\cdot)\eta(\cdot)$.
	\end{enumerate}   
	Given such a structure one can consider the almost complex structure $J$ on the Riemannian cone $\big( M\times\R^+,t^2g+\di t^2\big)$ given by
	\begin{enumerate}
		\item[$\bullet$] $J=\Phi$ on $\D=\ker\eta$, and
		\item[$\bullet$] $R=J(t\partial_t)_{\vert_{\{t=1\}}}$.
	\end{enumerate} 
	A \textit{Sasakian} \textit{structure} is a K-contact structure $(\eta,\Phi,R,g)$ such that the associated almost complex structure $J$ is integrable, and therefore $\left( M\times\R^+,t^2g+\di t^2,J\right)$ is K\"ahler. 
	A \textit{Sasakian} manifold is a manifold $M$  equipped with a Sasakian structure $(\eta,\Phi,R,g)$.
	
	Equivalently, one can define Sasakian manifolds in terms of \K cones. 
	Namely, a Sasakian structure on a smooth manifold $M$ is defined to be a \K cone structure on $M\times\R^+=Y$, that is, a \K structure $(g_Y,J)$ on $Y$ of the form 
	$g_Y=t^2g+\di t^2$ where $t$ is the coordinate on $\R^+$ and $g$ a metric on $M$.
	Then $(Y,g_Y,J)$ is called the \K cone of $M$ which, in turn, is identified with the submanifold $\{t=1\}$. The \K form on $Y$ is then given by
	$$\Omega_Y=\dfrac{i}{2}\de\deb t^2.$$
	The Reeb vector field on $Y$ is defined as 
	$$R=J(t\de_t).$$
	This defines a holomorphic Killing vector field with metric dual $1$-form
	$$\eta=\dfrac{g_Y(R,\cdot)}{t^2}=\di^c\log t=i(\deb-\de)\log t$$
	where $d^c=i(\deb-\de)$.
	Notice that $J$ induces an endomorphism $\Phi$ of $\TM$ by setting 
	\begin{enumerate}
		\item[$\bullet$] $\Phi=J$ on $\D=\ker\eta_{\vert_{\TM}}$, and
		\item[$\bullet$] $\Phi(R)=0$.
	\end{enumerate} 
	Equivalently, the endomorphism $\Phi$ is determined by $g$ and $\eta$ by setting
	$$g(X,Z)=\dfrac{1}{2}\di\eta(X,\Phi Z)\ \ \mbox{ for }\ X,Z\in\D.$$
	It is easy to see that, when restricted to $M=\{t=1\}$, $(\eta,\Phi,R,g)$ is a Sasakian in the contact metric sense whose \K cone is $(Y,g_Y,J)$ itself. 
	When this does not lead to confusion, we will use $R$ and $\eta$ to indicate both the objects on $Y$ and on $M$.

	Since $g$ and $\eta$ are invariant for $R$, the Reeb foliation is transversally \K in the sense that the distribution $\D$ admits a \K structure $(g^T,\omega^T,J^T)$ which is invariant under $R$.
	Explicitly, we have 
	$$\omega^T=\dfrac{1}{2}\di\eta,\ \ J^T=\Phi_{\vert_{\D}}\ \mbox{ and } \ g^T(X,Z)=\dfrac{1}{2}\di\eta(X,J^T Z)=g_{\vert_{\D}}.$$
	In particular, one can see that
	\begin{equation}\label{EqTransvFormIsCurv}
		\omega^T=\dfrac{1}{2}\di\eta=\dfrac{i}{2}\di(\deb-\de)\log t=i\de\deb\log t.
	\end{equation}
	
	The Reeb vector field defines a foliation on $M$, called the Reeb foliation.
	A very important dichotomy of Sasakian structures is given by the regularity of the leaves of the Reeb foliation.
	Namely, if there exist foliated charts such that each leaf intersects a chart finitely many times, the Sasakian structure is called \textit{quasi-regular}.
	Other wise it is called \textit{irregular}. Moreover, if every leaf intesects every chart only once, the sasakian structure is said to be \textit{regular}.
	Compact regular and quasi-regular Sasakian manifolds are fairly well understood due to the following structure theorem.
	
	\begin{theor}[\cite{Boyer08Book}]\label{TheorStructure}
		Let $(M,\eta,\Phi,R,g)$ be a quasi-regular  compact Sasakian manifold. Then the space of leaves of the Reeb foliation $(X,\omega,g_\omega)$ is a compact \K cyclic orbifold with integral \K form $\frac{1}{2\pi}\omega$ so that the projection $\pi:(M,g)\lra(X,g_\omega)$ is a Riemannian submersion. 
		Moreover, $X$ is a smooth manifold if and only if the Sasakian structure on $M$ is regular. 
		
		Viceversa, any principal $S^1$-orbibundle $M$ with Euler class $-\frac{1}{2\pi}[\omega]\in H^2_{orb}(X,\Z)$ over a compact \K cyclic orbifold $(X,\omega)$ admits a Sasakian structure.
	\end{theor}
	This result allows us to reformulate the geometry of a compact quasi-regular Sasakian manifold $M$ in terms of the  algebraic geometry of the \K orbifold $X$. 	
	We will illustrate in detail this correspondence for its importance in the remainder of the paper. 
	We only present this in the regular setting since it is the only one needed for our results. 
	Nonetheless, the same description applies to quasi-regular compact Sasakian manifolds.
	Let us first introduce the concept of $\D$-homothetic transformation of a Sasakian structure.
	
	\begin{defin}[$\D$-homothety or a transverse homothety]\label{DefinTransverseHomothety}
		Let $(M,\eta,\Phi,R,g)$ be a (not necessarily compact) Sasakian manifold and $a\in\R$ a positive number. One can define the Sasakian structure $(\eta_a,\Phi_a, R_a,g_a)$ from $(\eta,\Phi,R,g)$  as	
		$$ \eta_a=a\eta,\ \ \ \Phi_a=\Phi,\ \ \ R_a=\dfrac{R}{a},\ \ \  g_a=ag+(a^2-a)\eta\otimes\eta=ag^T+\eta_a\otimes\eta_a.$$
		
		Equivalently, we can define the same structure on $M$ by setting a new coordinate on the \K cone as $\widetilde{t}=t^a$. It is clear from the formulation above that this induces on $M=\{\widetilde{t}=1\}=\{t=1\}$ the same Sasakian structure $(\eta_a,\Phi_a, R_a,g_a)$. We will call this structure the $\D_a$-homothety of $(\eta,\Phi,R,g)$.
	\end{defin}
	
	Now let the compact regular Sasakian manifold $(M,\eta,\Phi,R,g)$ be given and consider the projection $\pi:(M,g)\lra(X,\omega)$ given in Theorem~\ref{TheorStructure} above. 
	Notice that $\pi$ locally identifies the contact distribution $\D$ with the tangent space of $X$. Therefore,  up to $\D$-homothety, we have that $\pi^*(\omega)=\frac{1}{2}\di\eta$. 
	Moreover, the endomorphism $\Phi$ determines the complex structure on $X$ and $g$ induces the \K metric $g_\omega$, i.e. $g^T=\pi^*g_\omega$.
	
	In this case the class $\frac{1}{2\pi}[\omega]\in H^2(X,\Z)$ defines an ample line bundle $L$ over $X$.
	Moreover, the cone $Y=M\times\R^+$ is identified with the complement of the zero section in $L^{-1}=L^*$ in the following way. Let $h$ be a hermitian metric on $L$ such that 
	$$\omega=-i\de\deb\log h.$$
	Then its dual $h^{-1}$ on $L^{-1}$ defines the second coordinate of $(p,t)\in M\times\R^+=L^{-1}\setminus\{0\}$ by
	\begin{align}
		\nonumber t:L^{-1}\setminus\{0\} &\lra \R^+\\ 
		\label{EqMetricOnConeFromHermitian} (p,v)&\mapsto\vert v\vert_{h^{-1}_p}
	\end{align}
	where $v$ is a vector of $L^{-1}$ in the fibre over $p$.
	With this notation the \K form on the \K cone $\left( M\times\R^+,t^2g+\di t^2,I\right)$ is given by
	
	\begin{equation}
		\Omega=\dfrac{i}{2}\de\deb t^2.
	\end{equation}
	The Sasakian structure can be read from this data as
	\begin{equation}\label{EqFromHermitianToSasakian}
		\omega^T=-i\de\deb\log h,\ \ \ \ \eta=i(\deb-\de)\log t.
	\end{equation}

	Therefore, the choice of a hermitian metric $h$ on an ample line bundle $L$ over a compact \K manifold $X$ completely determines a Sasakian structure on the $U(1)$-orbibundle associated to $L^{-1}$. 
	The Sasakian manifold $(M,\eta,R,g,\Phi)$ so obtained is called a \textit{Boothby-Wang bundle} over $(X,\omega)$. 
	Observe that, although the differentiable manifold is uniquely determined by $\frac{1}{2\pi}[\omega]$, the Sasakian structure does depend on, and is in fact determined by, the choice of $h$.

	The most basic example is the standard Sasakian structure on $S^{2n+1}$, that is, the Boothby-Wang bundle determined by the Fubini-Study metric $h=h_{FS}$ on $\OO(1)$ over $\CP^n$. We give the details of this construction to further illustrate the formulation above.
	\begin{ex}[Standard Sasakian sphere]\label{ExStandardSphere}
		Let $h=h_{FS}$ be the Fubini-Study hermitian metric on the holomorphic line bundle $\OO(1)$ over $\CP^n$. Recall that its dual metric $h^{-1}$ on $\OO(-1)\setminus\{0\}=\C^{n+1}\setminus\{0\}$ is given by the euclidean norm.
		This defines a coordinate $t$ on the \K cone $\OO(-1)\setminus\{0\}=\C^{n+1}\setminus\{0\}=S^{2n+1}\times \R^+$. Namely, for coordinates $z=(z_0,z_1,\ldots,z_n)$ on $\C^{n+1}$ we have
		\begin{align*}
			t:\C^{n+1} &\lra \R^+\\ 
			z&\mapsto\vert z\vert=\sqrt{\sum_{i=0}^{n}z_i\overline{z}_i}
		\end{align*}
		Now the \K metric on the cone is nothing but the flat metric
		$$	\Omega_{flat}=\dfrac{i}{2}\de\deb t^2=\dfrac{i}{2}\sum\di z_i\wedge\di \overline{z}_i.$$
		The Reeb vector field $R_0$ and the contact form $\eta_0$ read
		$$	R_0=J(t\de_t)=i\sum z_i\de_{z_i}-\overline{z}_i\de_{\overline{z}_i},\ \ \ \eta_0=i(\deb-\de)\log t=\dfrac{i}{2t^2}\sum z_i\di\overline{z}_i-\overline{z}_i\di z_i.$$
		It is clear that, when restricted to $S^{2n+1}$, $\eta_0$ and $R_0$, together with the round metric $g_0$ and the restriction $\Phi_0$ of $J$ to $\ker\eta_0$ give a Sasakian structure on $S^{2n+1}$.
		This corresponds exactly to the Hopf bundle $S^{2n+1}\lra\CP^n$. Moreover, we have
		$$\pi^*\omega_{FS}=\omega^T=\dfrac{1}{2}\di\eta_0=\dfrac{i}{2\vert z\vert^4}\sum_i\vert z_i\vert^2\di z_i\wedge\di \overline{z}_i-\sum_{i,j}\overline{z}_i z_j\di z_i\wedge\di \overline{z}_j$$
		where $\pi:\C^{n+1}\setminus\{0\}\lra\CP^n$ is the standard projection.
	\end{ex}
	
	In general the space of leaves of the Reeb foliations $X$ is not even an orbifold. Nevertheless, when the Sasakian structure is regular and complete, $X$ is a \K manifold regardless of the compactness of $M$, see e.g. \cite{Reinhart59Foliated}.
	
	We now switch our attention back to not necessarily compact Sasakian manifolds and recall another well known class of deformations of Sasakian structures, the so-called transverse \K transformations.
	Namely, given a \K cone $Y=M\times\R^{+}$ we consider all \K metrics on $(Y,J)$ that are compatible with the Reeb field $R$. In other terms, these are potentials $\widetilde{t}^2$ such that $t\de_t=\widetilde{t}\de_{\widetilde{t}}$.
	This means that the corresponding \K and contact forms satisfy 
	$$\widetilde{\Omega}=\Omega+i\de\deb e^{2f},\ \ \ \widetilde{\eta}=\eta+d^cf$$
	for a function $f$ invariant under $\de_t$ and $R$. Such functions are colled \textit{basic functions}.
	We still need to identify the manifolds $\{\widetilde{t}=1\}$ and $\{t=1\}$.
	This is done by means of the diffeomorphism
	\begin{align*}
		F:Y&\lra Y\\
		(p,t)&\mapsto\left(p,te^{-f(p)}\right)
	\end{align*}
	which maps $\{t=1\}$ to $\{t=e^{-f(p)}\}=\{\widetilde{t}=1\}$. 
	It is elementary to check that $\eta,\ R$ and $d^cf$ are invariant under $F$ so that $\widetilde{\eta}=\eta+d^cf$ holds on $M$.
	Furthermore, the transverse \K forms are related by $\widetilde{\omega}^T=\omega^T+i\de\deb f$.
	Notice that when the Sasaki structure is quasi-regular basic functions correspond to function on the base orbifold $X$. 
	Thus, if $t$ comes from a hermitian metric $h^{-1}$ on $L^{-1}$, such a transformation is given by replacing $h^{-1}$ with $e^fh^{-1}$ for a function $f:X\lra\C$ such that $\omega+i\de\deb f>0$. 
	This is equivalent to picking a different \K form $\widetilde{\omega}$ in the same class as $\omega$.
	We summarise the above discussion in the following definition.
	\begin{defin}[Transverse \K deformations]\label{DefinTransverseDeformation}
		Let  $(M,\eta,R,g,\Phi)$ be a Sasakian manifold with \K cone $(Y,J)$ and \K potential $t^2$. 
		A transverse \K transformation is given by replacing $t$ with $\widetilde{t}=e^{f}t$ for a basic function $f$ and leaving $(Y,J,R)$ unchanged. 
		When the Sasaki structure is quasi-regular and given as in \eqref{EqFromHermitianToSasakian}, a transverse \K transformation is given by replacing $h^{-1}$ with $e^fh^{-1}$.
	\end{defin}
	
	In this paper we are mostly interested with immersions and embeddings of Sasakian manifold so we recall the most relevant definitions.
	Two Sasakian manifolds $(M_1,\eta_1,R_1,g_1,\phi_1)$ and  $(M_2,\eta_2,R_2,g_2,\phi_2)$ are \text{equivalent} if there exists a diffeomorphism $\varphi\colon M_1\lra M_2$ such that
	\begin{equation*}
		\varphi^*\eta_2=\eta_1\hspace{5mm}\mbox{and}\hspace{5mm}\varphi^*g_2=g_1.
	\end{equation*}
	If this holds then $\varphi$ also preserves the endomorphism $\phi_1$ and the Reeb vector field.
	As implicitly intended in the definitions above, a Sasakian equivalence from a Sasakian manifold $(M,\eta,R,g,\phi)$ to itself is often called a \textit{Sasakian transformation} of $(M,\eta,R,g,\phi)$.

	\subsection{Complex and Sasakian space forms}
	In order to fix the notation we recall some facts on complex and Sasakian space forms. In the literature a space form is not assumed to be simply connected and complete but for us this will always be the case. 
	\begin{defin}
		A \textit{complex space form} $F(N,b)$, $N\leq\infty$, is a complete simply connected complex manifold of real dimension $2N$ with constant holomorphic sectional curvature $4b$.
	\end{defin}
	
	Complex space forms are distinguished by the sign of the holomorphic sectional curvature. Namely, up to holomorphic isometry, there are three types of complex space forms.
	\vskip 1em
	\noindent		\textbf{($b=0$)} $F(N,0)$ is the complex space $\C^N$. If $N=\infty$, $F(\infty,0)$ is the Hilbert space $\ell^2(\C)$ of sequences $(z_1,z_2,\ldots)$ of complex numbers such that $\sum_{j=1}^\infty \vert z_j\vert^2<\infty$. In both cases $F(N,0)$ is considered to be endowed with the flat metric.
	\vskip 0.3em
	\noindent		\textbf{($b>0$)} $F(N,b)$ is the complex projective space $\CP^N$.
	When $N=\infty$, $F(\infty,b)=\CP^\infty$ is the space of equivalence classes $[z_1:z_2:\cdots]$ such that $(z_1,z_2,\ldots)\in\ell^2(\C)\setminus\{0\}$ with respect to standard $\C^*$-action. This space is endowed with the Fubini-Study form $\omega_{FS}$ of constant holomorphic sectional curvature $4b$, given in homogeneous coordinates by $$\omega_{FS}=\dfrac{i}{2b}\de\deb\log \sum_{j=1}^N \vert z_j\vert^2.$$
	\vskip 0.3em
	\noindent		\textbf{($b<0$)} $F(N,b)=\CH^\infty$ is the ball of radius $-\frac{1}{b}$ in $\C^N$ (or $\ell^2(\C)$ when $N=\infty$)  endowed with the hyperbolic \K form $\omega_{hyp}$ of constant holomorphic sectional curvature $4b$:
	$$\omega_{hyp}=\dfrac{i}{2b}\de\deb\log \left(1+b\sum_{j=1}^N \vert z_j\vert^2\right).$$
	\vskip 1em

	Sasakian space forms are defined analogously. Given a Sasakian manifold $(M,\eta,R,g,\phi)$ with Riemannian sectional curvature $\mathrm{Sec}$, the \textit{$\phi$-sectional curvature $H$} of $g$ is defined by 
	$$H(X)=\mathrm{Sec}(X,\phi X)$$
	for all vector fields $X$ of unit length orthogonal to $R$.
	\begin{defin}
		A \textit{Sasakian space form} $M(N,c)$, $N\leq\infty$, is a complete simply connected Sasakian manifold of dimension $2N+1$ with constant $\phi$-sectional curvature $H\equiv c$.
	\end{defin}

	Tanno \cite{Tanno69SpaceForm} proved that there are three types of finite dimensional Sasakian space forms, namely, those with $ c<-3$, $c=-3$ and $c>-3$.
	Under the Boothby-Wang correspondence, constant $\phi$-sectional curvature $c$ corresponds precisely to constant holomorphic transverse sectional curvature $c+3$. 
	Namely, the $\phi$-sectional curvature is related to to the sectional curvature $\mathrm{Sec}^T$ of the transverse metric $g^T$ by
	\begin{equation}\label{EqSecCurv}
		H(X)=\mathrm{Sec}^T(X,\phi X)-3
	\end{equation}
	for a vector field $X$ of unit length orthogonal to $R$, see for instance \cite[Section~7.3]{Boyer08Book}. 
	We report Tanno's classification including the infinite dimensional case with particular focus on the relation with complex space forms.
	\vskip 1em
	\noindent		\textbf{($c=-3$)} $M(N,-3)$ is Sasaki equivalent to $\R^{2N+1}=\R\times F(N,0)$ with the contact standard structure. In particular, if $N=\infty$, then $M(\infty,-3)$ is the Hilbert space $\R\times\ell^2(\C)$. In both cases $M(N,-3)$ is equipped with the contact form $\eta=\di t+\frac{i}{2}\sum_{j=1}^N (z_j\di\overline{z}_j-\overline{z}_j\di z_j)$ and Reeb vector field $\de_t$ where $t$ is the coordinate on the first factor. The endomorphism $\phi$ is defined to be the horizontal lift of the complex structure of $\C^N$ ($\ell^2(\C)$) via the projection $\pi:\R\times F(N,0)\lra F(N,0)$. The metric $g$ is defined to be $g=g^T+\eta\otimes\eta$ where $g^T$ is the flat metric on $F(N,0)$.
	\vskip 0.3em
	\noindent		\textbf{($c>-3$)} $M(N,c)$ is Sasaki equivalent to the standard Sasakian sphere $S^{2N+1}$ ($S^{\infty}$ if $N=\infty$). This is the Boothby-Wang bundle over $\CP^N=F(N,b)$ whose transverse \K metric has constant holomorphic sectional curvature $4b=c+3$.
	\vskip 0.3em
	\noindent		\textbf{($c<-3$)} $M(N,c)$ is Sasaki equivalent to $\R\times F(N,b)$ where the transverse K\"ahler structure is that of the hyperbolic complex space $F(N,b)$ of constant holomorphic sectional curvature $4b=c+3$. Namely, consider the projection on the second factor $\pi:\R\times F(N,b)\lra F(N,b)$ and let $\alpha$ be a $1$-form on $F(N,b)$ such that $\di\alpha=\omega_{hyp}$.
	Define the contact form of $M(\infty,c)$ to be $\eta=\di t+\pi^*\alpha$ with Reeb vector field given by $\de_t$ where $t$ is the coordinate on the first factor. The endomorphism $\phi$ is defined to be the horizontal lift of the complex structure of $F(N,b)$. The metric $g$ is given by $g=g^T+\eta\otimes\eta$ where $g^T$ is the Bergman metric on $F(N,b)$ associated to $\omega_{hyp}$.
	\vskip 1em

	The choice of the curvature being $4b$ instead of $b$ derives from the fact that, with the normalisation we have fixed, the Fubini-Study form on $\CP^N$ has constant sectional holomorphic curvature $4$.
	Moreover, with this choice the sphere $S^{2N+1}$ with the round metric is the Sasakian space form $M(N,1)$ fibring over the complex space form $\CP^N=F(N,1)$.
	
	Notice that each simply connected space form $M(N,c)$ admits a fibration $\pi:M(N,c)\rightarrow F(N,b)$ over a complex space form with $4b=c+3$ whose fibres are the leaves of the Reeb foliation. When $c>-3$,  $\pi$ is a principal circle bundle, while in the other two cases the fibration is trivial.
	Observe that in these two cases, any choice of any primitive $\alpha$ of the \K form of the complex space form $F(N,b)$ yields a Sasakian structure of constant $\phi$-sectional curvature which is Sasaki equivalent to the Sasakian structure we have fixed. 
	
	\subsection{Some remarks on noncompact Sasakian manifolds}
	We conclude this section by discussing some basic properties of noncompact regular Sasakian manifolds and, in particular, of noncompact homogeneous Sasakian manifolds.
	These are fibre bundles over homogeneous \K manifolds with fibre $\R$ or $S^1$, see \cite[Theorem~8.3.5]{Boyer08Book}. 
	Let $M$ be a homogeneous Sasakian manifolds and $X$ be the \K base of the associated fibre bundle. 
	
	Firstly, we point out that the bundle is trivial if the fibre is $\R$. This follows from the fact that both the total space $M$ and the base $X$ are orientable, therefore the bundle is too. Since this is real line bundle, orientable means trivial.

	Now we want to show that if the base manifold $X$ is compact, then the Sasakian manifold $M$ cannot be the total space of a trivial bundle. 
	This follows directly from the structure theorem if $M$ is compact. 
	For $M$ noncompact it follows from the next
	
	\begin{lem}\label{LemNonCompact}
		Let $X$ be a \K manifold. The product $X\times \R$ admits a Sasakian structure such that the Reeb orbits are the fibres of the trivial fibration if and only if the transverse \K form $\omega$ is exact on $X$.
	\end{lem}
	\begin{proof}
		We write $M=X\times\R$ and $\pi:X\times\R\lra X$ for the projection onto the first factor. Let $t$ be the coordinate on the second factor.
		
		Suppose firstly that $\omega=\di\alpha$ is an exact \K form on $X$ and $g^T$ is the associated metric. Moreover, let $t$ be the coordinate on the second factor.
		Then the form $\eta=\di t+\pi^*\alpha$ is a contact form on $M$ defining a Sasakian structure $(\eta,\de_t,g,\Phi)$ with $g=g^T+\eta\otimes\eta$ and $\Phi$ the horizontal lift of the complex structure on $X$. 
		In this case the projection onto the transverse \K space is given simply by $\pi$ above.
		
		Suppose now that $M$ admits a Sasakian structure $(\eta,R,g,\Phi)$ such that the Reeb orbits are the fibres of the trivial fibration. That is, the Reeb orbit through $(x,t)$ is $\{x\}\times\R$ and $\pi$ above is exactly the projection onto the space of leaves of the Reeb foliation.
		Then the flow of the Reeb vector field defines a coordinate $t:\R\lra X\times\R$ for the second factor so that $R=\de_t$.
		Then we can write $\eta=\di t+\gamma$ where $\gamma(\de_t)=0$. 
		Moreover, since $\iota_{\de_t}\di\gamma=\iota_{\de_t}\di\eta=0$, the form $\gamma$ is basic, that is, $\gamma=\pi^*\alpha$ for some $1$-form $\alpha$ on $X$. Therefore $\pi^*\omega=\di\eta=\di\gamma=\di\pi^*\alpha$ is exact, and so is $\omega$.
	\end{proof}
	This result also holds if we replace $\R$ by $S^1$. 
	In fact, such a Sasakian structure on $X\times S^1$ can be lifted to $X\times \R$. 
	Viceversa, if $X$ admits a \K form $\omega=\di\alpha$, then, for any volume form $\beta$ on $S^1$, the contact form $\eta=\beta+\pi^*\alpha$ defines a Sasakian structure on $X\times S^1$.
	In particular, if $X$ is compact, then Lemma~\ref{LemNonCompact} implies that there cannot exist any Sasakian structure on $X\times\R$ fibring over $X$. Clearly, since $X\times S^1$ is compact, this is also given by Theorem~\ref{TheorStructure}.
	For future reference, we specialise the discussion above to homogeneous manifolds.
	\begin{prop}\label{PropSasakiNonCompactfibreR}
		If $M$ be a noncompact homogeneous Sasakian manifold, then $M$ is the total space of a fibration over a noncompact homogeneous \K manifold $X$. 
		Moreover, if the fibre is $\R$, then $X$ is contractible.
	\end{prop}
	\begin{proof}
		Recall that a homogeneous Sasakian manifold $M$ is regular and fibres over a homogeneous \K manifold $X$ with fibre $\R$ or $S^1$. 
		If the base $X$ is compact, then the fibre is forced to be $\R$ because $M$ is noncompact. 
		By orientability, the fibration is trivial and Lemma~\ref{LemNonCompact} implies that $X$ admits an exact \K form.
		We conclude that the base $X$ must be noncompact. 
		
		If, additionally, the fibre is $\R$, once again triviality of the bundle follows from the orientability of $M$ and $X$. 
		Then Lemma~\ref{LemNonCompact} implies again that $X$ admits an exact \K form.
		Therefore $X$ admits a global \K potential, that is, $X$ is a contractible manifold, see \cite[Theorem~1.2]{Loi15HomogeneousImmersions}.
	\end{proof}
	
	Notice that the assumption on the fibre of the fibration is necessary to get the contractability of $X$.
	For instance, one can construct a homogeneous Sasakian structure on the circle bundle over the product $\CP^k\times \CH^l$.
	This argument shows that any homogeneous Sasakian manifold fibring over a homogeneous \K manifold $X$ without compact factors is $X\times\R$ with the structure described above.

	\section{Extensions  and uniqueness of Sasakian immersions in Sasakian space forms}\label{SecRigidity}
	In this section we give a general discussion on immersions of Sasakian manifolds into Sasakian space forms. Given two Sasakian manifolds $(M_1,\eta_1,R_1,g_1,\phi_1)$ and  $(M_2,\eta_2,R_2,g_2,\phi_2)$, a \textit{Sasakian immersion (embedding)} of $M_1$ in $M_2$ is an immersion (embedding) $\varphi\colon M_1\lra M_2$ such that
	\begin{align*}
		\varphi ^*\eta_2&=\eta_1,\hspace{20mm}\varphi ^*g_2=g_1,\\
		\varphi_*R_1&=R_2\hspace{5mm}\mbox{and}\hspace{5mm}\varphi_*\circ \phi_1=\phi_2\circ\varphi_*.
	\end{align*}  
	We can rephrase this definition in terms of the \K cone of the Sasakian manifolds $M_1$ and $M_2$. Namely, the map $\varphi$ satisfying the conditions above clearly extends to a map
	\begin{align*}
		\widetilde \varphi\colon M_1\times\R&\lra M_2\times\R\\
		(p,t)&\mapsto (\varphi(p),t).
	\end{align*}  
	It is clear that if $\varphi$ is a Sasakian immersion (resp. embedding), then $\widetilde \varphi$ is a \K immersion (resp. embedding).
	
	If, conversely, $Y_1$ and $Y_2$ are the \K cones of $M_1$ and $M_2$ with coordinates $t_1$ and $t_2$, then a \K immersion (embedding) $\widetilde \varphi\colon Y_1 \lra Y_2$ such that $\widetilde\varphi^*(t_2)=t_1$ restricts to a Sasakian immersion (embedding) $\varphi:M_1\lra M_2$.
	Since it is often more useful to our purposes, we give the following 
	
	\begin{defin}[Sasakian immersion and embedding]\label{DefEmbedding}
		Let $M_1$ and $M_2$ be two Sasakian manifolds with \K cones $Y_1$ and $Y_2$ and coordinates $t_1$ and $t_2$ respectively. 
		A \textit{Sasakian immersion (respectively embedding)} of $M_1$ in $M_2$ is a \K immersion (resp. embedding) $\varphi\colon Y_1\lra Y_2$ such that $\varphi^*(t_2)=t_1$.
	\end{defin}
	\begin{remark}\rm
		Given the equivalence between a Sasakian immersion $M_1\lra M_2$ and a \K immersion of the \K cones,  with an abuse of notation, we will often denote both maps with the same letter. 
	\end{remark}
	
	\begin{proof}[Proof of~\Cref{ThmSasakianRigidity}] We will prove the theorem in the case where $M$ is a trivial fibration over contractible \K manifold $X$. We then apply this to Darboux neighbourhoods in order to get the thesis by a standard argument.
		
		Let $M$ be a Sasakian manifold such that the Reeb foliation defines a map to a contractible \K manifold $X$, that is, $M$ is diffeomorphic to $X\times\R$.
		Firstly, we focus on the case where $c>-3$. Here the \K cone $M\times\R^+$ admits two \K immersions $\widetilde\varphi_1=\varphi_1\times \Id_{\R^+}$ and $\widetilde\varphi_2=\varphi_2\times \Id_{\R^+}$ into $S^{2N+2}\times\R^+=\C^{N+1}\setminus\{0\}$ and therefore into the complex space form $\C^{N+1}$. 
		Hence, by Calabi's rigidity, there exists a rigid transformation $F$ of $\C^{N+1}$ such that $\widetilde\varphi_1=F\circ\widetilde\varphi_2$.
		This yields the thesis as $F$ is forced to be unitary by the very definition of $\widetilde \varphi_1$ and $\widetilde \varphi_2$.

		We can now assume $c\leq-3$. 
		Recall that in this case  $M(N,c)=F(N,b)\times\R$ with $4b=c+3$.
		Since Sasakian immersions map Reeb leaves to Reeb leaves, the immersions $\varphi_1,\varphi_2$ cover two \K immersions defined so that the following diagram commutes
		$$
		\begin{tikzcd}[column sep=large, row sep=large]
			M=X\times\R \arrow [r,"\varphi_i"] \arrow[d, "\pi"] &
			M(N,c) \arrow[d,"\widetilde \pi"] \\
			X \arrow[r,"\psi_i"] & F(N,b)
		\end{tikzcd}
		$$
		for $i=1,2$.
		Since the fibrations $\pi$ and $\widetilde \pi$ are trivial, the maps $\varphi_i$ have the form $\varphi_i(x,t)=(\psi_i(x),f_i(x,t))$ for $x\in X$, $f_i\in C^\infty(M,\R)$ and $i=1,2$.
		Now, by Calabi's rigidity, there exists a rigid transformation $\hat T$  of $F(N,b)$ such that $\psi_1=\hat T\circ\psi_2$. 
		This lifts to a rigid transformation $\widetilde T :M(N,c)\lra M(N,c)$ given by $\widetilde T(p,s)=(\hat T (p),s)$.
		Therefore we get an immersion $\hat \varphi_2:=\widetilde T\circ\varphi_2:M\lra M(N,c)$ defined as $\hat \varphi_2(x,t)=(\hat T\circ\psi_2(x),f_2(x,t))=(\psi_1(x),f_2(x,t))$.
		
		Now, the fact that the immersions $\varphi_1$ and $\varphi_2$ are Sasakian yields  $f_i(x,t)=t+a_i$ for $i=1,2$. 
		Let $a=a_1-a_2$ and denote by $T_a:M(N,c)\lra M(N,c)$ the translation $T_a(p,s)=(p,s+a)$.
		Then the Sasakian immersions $\varphi_1$ and $\varphi_2$ are related by $\varphi_1= T\circ\varphi_2$ where $T$ is the rigid transformation given by $T=T_a\circ\widetilde T$. 
		
		Notice that the transformation $T$ is unique when the immersion is full, that is, when its image is not contained in a totally geodesic submanifold of $M(N,c)$. This follows from the uniqueness of the rigid transformation $\hat T$, see the proof of~\cite[Theorem~2.1.4]{Loi18Book}.
		Now let $M$ be an arbitrary Sasakian manifold and $\varphi_1,\varphi_2$ two immersions of $M$ into a Sasakian space form $M(N,c)$.
		By restricting the target to a totally geodesic submanifold of $M(N,c)$, we can assume that $\varphi_1,\varphi_2$ are full.
		Each point $p\in M$ has a Darboux neighbourhood $U_p$ that can be regarded as a trivial Sasakian fibration over a open disc in $\C^n$ (where $\dim(M)=2n+1$). 
		The restrictions of $\varphi_1,\varphi_2$ to $U_p$ define Sasakian immersions of $U_p$ into $M(N,c)$ which, by the first part of this proof, satisfy $\varphi_{1\vert U_p}= T_p\circ\varphi_{2\vert U_p}$ for a unique rigid transformation $T_p$ of $M(N,c)$.
		This defines a map $p\mapsto T_p$ from $M$ to the group of rigid transformations of $M(N,c)$.
		This map is constant because it is locally constant and $M$ is connected.
		Therefore there exists a unique Sasakian transformation $T$ of $M(N,c)$ such that $\varphi_1= T\circ\varphi_2$ globally.
	\end{proof}
	
	Now, using~\Cref{ThmSasakianRigidity}, we are able to prove~\Cref{ThmSasakianGluing} by a classical argument

	\begin{proof}[Proof of~\Cref{ThmSasakianGluing}]
		Consider, for each $p\in M$, a geodesically convex Darboux chart $W_p\subset U_p$ such that the collection $\{W_p\}_{p\in M}$ is still an open covering of $M$.
		Now fix a point $p_0\in M$. For any point $p\in M$ there exists a path $\gamma$ from $p_0$ to $p$. 
		By compactness there exist a finite number of charts $W_{p_0},\ldots,W_{p_k}$ that cover the image of $\gamma$. 
		Moreover, we can order the $W_{p_i}$'s so that the intersections $W_{p_i}\cap W_{p_{i+1}}\cap\gamma$ are nonempty for all $i=0,\ldots,k-1$.
		Denote by $\gamma_i$ the intersection of $W_{p_i}$ with $\gamma$.
		Now the set $W_{p_0}\cap W_{p_1}$ admits two Sasakian immersions into $M(N,c)$.
		By Theorem~\ref{ThmSasakianRigidity} there exists a unique Sasakian transformation $T_1$ of $M(N,c)$ such that $T_1\circ \varphi_{p_1}=\varphi_{p_0}$ when restricted to $W_{p_0}\cap W_{p_1}$.
		Thus, by taking $T_1\circ \varphi_{p_1}$ in place of $\varphi_{p_1}$, we can extend the Sasakian immersion to $\gamma_0\cup \gamma_1$.
		By repeating this process a finite number of times we extend the immersion to the entire path $\gamma$ and, since $M$ is path connected, to the whole manifold $M$.
		Moreover, this extension is well defined because $M$ is simply connected.
		The last statement follows directly from Theorem~\ref{ThmSasakianRigidity}.
	\end{proof}

	\section{Sasakian immersions of homogeneous Sasakian manifolds into Sasakian space forms with $c\leq-3$}\label{SecImmersionNegative}
	
	In this section we investigate immersions into Sasakian space forms in the special case of homogeneous or, more generally, locally homogeneous Sasakian manifolds.
	This is inspired by the parallel with \K geometry.
	In particular, we seek a Sasakian analogue of the following result.
	\begin{theor}[\cite{Discala12Immersions}]\label{mainHomKahler}
		Let $(X,\omega)$ be a complex $n$-dimensional homogeneous K\"ahler manifold and suppose that $X$ admits a \K embedding into $F(N,b)$ with $N\leq\infty$.
		\begin{enumerate}
			\item If $b<0$ then $X=F(n,b)$. \label{KahlerImmersionHyperbolic}
			\item If $b=0$ then $X=F(k,0)\times F(n_1,b_1)\times\cdots\times F(n_r,b_r)$ where $k+n_1+\dots+n_r=n$ and $b_i<0$ for all $i=1,\ldots,r$.
		\end{enumerate}
	\end{theor}
	Theorem~\ref{mainHomKahler} will be crucial in the following proof of its Sasakian analogue.
	\begin{proof}[Proof of~\Cref{mainHom}]
		\noindent \eqref{point1mainHom} Assume $c<-3$ and consider the universal cover $\widetilde \pi:\widetilde M\lra M$ of $M$. Since $M$ is complete and locally homogeneous its universal cover $\widetilde M$ is a homogeneous Sasakian manifold.
		Notice that $\widetilde M$ admits a local Sasakian immersion in $M(N, c)$ at any point in the preimage $\widetilde\pi^{-1}(U)$.
		Therefore, $\widetilde M$ admits an embedding at any point because its group of Sasakian transformation acts transitively.
		Now, by Theorem~\ref{ThmSasakianGluing}, there exists a global Sasakian immersion $\widetilde\varphi:\widetilde M\lra M(N, c)$.

		Recall that $\widetilde M$ is regular and it is a fibration over a homogeneous K\"ahler manifold $(\widetilde X,\omega)$ whose fibre is either $\R$ or $S^1$, cf. \cite[Theorem~8.3.5]{Boyer08Book}. 
		Therefore $\widetilde \varphi$ covers a \K immersion $\widetilde \psi:\widetilde X\lra F(N,b)$ with $4b=c+3$.
		By Theorem~\ref{mainHomKahler} this implies that $\widetilde{X}$ is biholomorphically isometric to the complex space form $F(n,b)$ and therefore the fibre of the bundle $\widetilde{M}\lra\widetilde{X}$ must be $\R$ because $\widetilde{X}$ is simply connected.
		From the fact that $\widetilde\varphi$ is Sasakian we conclude that the bundle $\widetilde \pi:\widetilde M\lra M$ is nothing but the pullback under $\widetilde\psi$ of the bundle $M(N,c)\lra F(N,b)$.
		In particular, $\widetilde{M}$ is Sasaki equivalent to $M(n,c)$.
		Namely, we get the following commutative diagram where horizontal maps are \K or Sasakian immersions and vertical maps are the standard fibrations
		$$
		\begin{tikzcd}[column sep=large, row sep=large]
			\widetilde{M}=M(n,c) \arrow [r,"\widetilde{\varphi}"] \arrow[d] &
			M(N,c) \arrow[d] \\
			\widetilde X =F(n,b) \arrow[r,"\widetilde{\psi}"] & F(N,b)
		\end{tikzcd}
		$$
		We conclude that $M$ is Sasaki equivalent to $M(n,c)/\Gamma$ where $\Gamma=\pi_1(M)$ acts by Sasakian transformations.

		\noindent \eqref{point2mainHom} The line of argument in the case where $c=-3$ is similar to the one followed in the previous point so we only point out the differences.
		
		In this case, by Theorem~\ref{mainHomKahler}, the homogeneous \K manifold $\widetilde{X}$ is holomorphically isometric to
		$F(k,0)\times F(n_1,b_1)\times\cdots\times F(n_r,b_r)$
		where $k+n_1+\dots+n_r=n$ and $b_i<0$ for all $i=1,\ldots,r$. 
		Therefore the manifold $\widetilde{M}$ is Sasaki equivalent to the standard Sasakian manifold over it. 
		In other words, $\widetilde{M}$ is the pullback under the immersion $\widetilde{\psi}$ of the bundle $M(N,c)\lra F(N,b)$. 
		One can also regard $\widetilde{M}$ as the Sasakian manifold obtained by identifying the fibres of the standard Sasakian manifolds on the factors of $\widetilde X$, that is, the Sasakian space forms $M(k,-3)$ and $M(n_i,4b_i-3)$ for $i=1,\dots, r$, which yields the thesis.
	\end{proof}

	\begin{remark}\label{RmkGlobalImmersion}\rm
		In the case where $U=M$, the immersion is in fact very explicit.
		This is because $M$ is Sasaki equivalent to $M(n,c)$ and the immersion $\widetilde\varphi$ covers the \K immersion $\widetilde\psi:\widetilde X\lra F(N,b)$, cf \cite[Theorem~1-2]{Loi18Book}.
		Namely, up to a rigid transformation of $M(N,c)$, when $c<-3$ the immersion $\widetilde\varphi$ is given by the linear inclusion of $M(n,c)=\R\times\CH^n$ into $M(N,c)=\R\times\CH^N$.
		When $c=-3$ we have $M=\R\times F(k,0)\times F(n_1,b_1)\times\cdots\times F(n_r,b_r)$. 
		In this case the immersion is given by $(f_0,f_1,\ldots,f_r)$ where $f_0$ is the linear inclusion of the first two factors and, for all $i=1,\ldots,r$,
		$$f_i(z_1,\ldots,z_{n_i})=\dfrac{-1}{2\sqrt{b_i}}\left(\ldots,\sqrt{\dfrac{(\vert j\vert-1)!}{j!}}z_1^{j_1}\cdots z_{n_i}^{j_{n_i}},\ldots\right)$$
		with $\vert j\vert=j_1+\cdots+j_{n_i}$ and $j!=j_1!\cdots j_{n_i}!$.
		In particular this implies that, if $N<\infty$ and $U=M$ in part \eqref{point2mainHom} of Theorem~\ref{mainHom}, then only the first two factors appear, i.e., $M=M(n,-3)$.
	\end{remark}
	\begin{remark}\label{RmkFundGroup}\rm
		By \cite[Theorem~8.3.5]{Boyer08Book} $M$ as in Theorem~\ref{mainHom} is regular and it is a fibration over a homogeneous K\"ahler manifold $X$ whose fibre is either $\R$ or $S^1$.
		The universal cover of $X$ is exactly the manifold $\widetilde{X}$ in the proof of Theorem~\ref{mainHom}.
		When the fibre of the bundle $M\lra X$ is $\R$, the group $\Gamma$ is a group of transverse \K transformations as it is exactly the fundamental group of $X=\widetilde{X}/\Gamma$.
		If instead the fibre is $S^1$, then $\Gamma$ is a central $\Z$-extension of the group $\pi_1(X)$ of \K isometries of $\widetilde{X}$.
		In particular, when $c=-3$, we have $$X= \dfrac{F(k,0)}{\Delta}\times \cfrac{F(n_1,b_1)}{\Gamma_1}\times\cdots\times\dfrac{F(n_r,b_r)}{\Gamma_r}$$ 
		where $\Delta$ (resp. $\Gamma_i$) is a discrete group of K\"ahler transformations of $F(k,0)$ (resp. $F(n_i,b_i)$ for all $i$). Namely, $\pi_1(X)$ splits as the direct product $\Delta\times\Gamma_1\times\cdots\times\Gamma_r$.
	\end{remark}

	\section{CR immersions of regular and complete Sasakian manifolds into spheres}\label{SecGeneralSasakianEmbedding}
	In the remainder of the paper we will study sufficient conditions for the existence of Sasakian embeddings or immersions of homogeneous Sasakian manifolds into the standard Sasakian sphere.
	In all the cases we consider, such embeddings are obtained via the same general method. 
	Namely, we use the algebraic geometry of the space of leaves of the Reeb foliation, that is, a \K manifold, to produce an embedding into a standard Sasakian sphere. 
	This section is devoted to the presentation of such CR embeddings for compact and noncompact regular Sasakian manifolds which are, in all respects, the Sasakian analogues of the Kodaira embedding. 
	Notice that CR embeddings of Sasakian manifolds into spheres have been studied also in \cite{Ornea07CREmbedding}. Nevertheless, here we discuss an explicit expression of the induced structure.
	
	Let $M$ be a compact regular Sasakian manifold. 
	By the Structure Theorem, $M$ is a $U(1)$-principal bundle $\pi:M\lra X$ over a compact \K manifold $(X,\omega)$ with $2\pi^*\omega=\di\eta$. 
	Furthermore, $M$ is the unitary bundle associated to the line bundle $L^{-1}$ where $c_1(L)=[\omega]$.
	This last condition implies that $L$ is ample. In other terms, $(X,L)$ is a polarised \K manifold.
	Therefore, for $k\in\N$ large enough, the bundle $L^{\otimes k}=L^k$ is very ample, and we can define the Kodaira embedding $\psi_k:X\lra \CP^{N_k}$ where $\dim(H^0(L))=N_k+1$.
	
	We want to construct a CR embedding $\varphi_k: M\lra S^{2N_k+1}$ of $M$ into the standard sphere covering the Kodaira embedding $\psi_k$. 
	Equivalently, by Definition~\ref{DefEmbedding} we can construct a holomorphic embedding of $\varphi_k:Y\lra \C^{N_k+1}\setminus \{0\}$ of the \K cone $Y=M\times\R^+$ into the \K cone $S^{2N_k+1}\times\R^+$.
	
	Let $h$ be a hermitian metric on $L$ whose Ricci curvature form is $\omega$ and such that the dual metric $h^{-1}$ on $L^{-1}$ induces the coordinate on $M\times \R^+=L^{-1}\setminus\{0\}$ as in \eqref{EqMetricOnConeFromHermitian}.
	This induces metrics $h^k$ and $h^{-k}$ on the line bundles $L^k$ and $L^{-k}=(L^*)^k$ respectively.
	Choose an orthonormal basis $(s_0, \ldots, s_{N_k})$ for the space of sections $H^0(L^k)$ with respect to the scalar product $\langle\cdot,\cdot\rangle_k$ defined by
	\begin{equation}\label{EqMetricOnSections}
		\langle s,t\rangle_k=\int_Xh^k\left(s(x),t(x)\right)\dfrac{\omega^n}{n!}
	\end{equation}
	where $s,t\in H^0(L^k)$.
	Fix an open set $U\subset X$ and a local trivialisation $\sigma:U\lra L^k$ so that we have $s_j=f_j\sigma $ for a holomorphic function $f_j:U\lra\C$.
	Then the Kodaira embedding is given locally by
	\begin{align*}
		\psi_{k\vert_U}: U&\lra\ \ \ \ \ \CP^{N_k}\\
		x&\mapsto \left[ \dfrac{s_0(x)}{\sigma(x)}:\cdots:\dfrac{s_{N_k}(x)}{\sigma(x)}\right].
	\end{align*}
	We have $\psi_k^*\OO(1)=L^k$, that is, $L^{-k}$ is the restriction of the tautological bundle $\OO(-1)$ to $\psi_k(X)$.
	Now $\sigma$ defines a trivialisation of $L^{-k}$ over $U$ so that a vector $v\in L^{-k}_{x}$ for $x\in U$ can be written as $\alpha_v=v(\sigma(x))$ with $\alpha_v\in\C$.
	Denote by $\widetilde\psi_k$ the holomorphic embedding $L^{-k}\setminus\{0\}\lra \C^{N_k+1}\setminus\{0\}$ induced by the Kodaira embedding.
	From the identities $\alpha_v=v(\sigma)$ and $ s_j=\sigma f_j$ yield the following local expression of $\widetilde\psi_k$.
	\begin{align*}
		U\times \C^*&\lra\ \ \ \ \ \ \ \ \ \C^{N_k+1}\setminus\{0\}\\
		(x,\alpha_v)&\mapsto (\alpha_v f_0(x), \ldots,\alpha_v  f_N(x))=\left( v( s_0(x)),\ldots,v( s_N(x))\right).
	\end{align*}

	Recall that taking the $k$-th power defines a holomorphic $k$-fold covering $p_k:L^{-1}\setminus\{0\}\lra L^{-k}\setminus\{0\}$.
	Thus we can define a holomorphic immersion 
	\begin{equation}\label{EqCompactEmbedding}
		\varphi_k=\widetilde\psi_k\circ p_k:L^{-1}\setminus \{0\}\lra \C^{N_k+1}\setminus \{0\}=\OO(-1)\setminus \{0\}.
	\end{equation}
	The coordinate free description of this embedding reads
	\begin{align*}
		\varphi_k:L^{-1}\lra&\ \ H^0(L)^*\\
		v\ \  \mapsto &\left(s\mapsto \dfrac{v(s)}{\vert s\vert}\right)
	\end{align*}
	where the norm $\vert s\vert$ is the one induced by \eqref{EqMetricOnSections}.
	Moreover, the pullback via $\varphi_k$ of the Fubini-Study hermitian metric $h^{-1}_{FS}$ on $\OO(-1)$ is computed as
	\begin{align*}
		\varphi_k^*(h^{-1}_{FS})_x(w,w) &=\sum_{j=0}^{N_k}\left\vert\alpha_{w^k}f_j(x)\right\vert^2=\sum_{j=0}^{N_k}\left\vert w^k(s_j(x))\right\vert^2=\sum_{j=0}^{N_k} h^{-1}(w^k,w^k)h^k\left(s_j(x),s_j(x)\right)\\
		&=h^{-k}(w,w)\sum_{j=0}^{N_k} h^k\left(s_j(x),s_j(x)\right)=B_k(x)h^{-k}(w,w)
	\end{align*}
	where $B_k(x)=\sum_{j=0}^{N_k} h^k\left(s_j(x),s_j(x)\right)$ is the Bergman kernel of $(X,L^k)$, $w\in L^{-1}\setminus\{0\}$ and $p_k(w)=w^k\in L^{-k}$.
	We summarise the discussion above in the following
	\begin{prop}\label{PropCompactEmbedding}
		Let $M$ be the compact regular Sasakian manifold determined by the Hermitian bundle $(L,h)$ over a compact projective manifold $X$. 
		For every integer $k>>0$ there exists a holomorphic embedding $\varphi_k:M\times\R^+\lra S^{2N_k+1}\times\R^+$ such that
		$\varphi_k^*(\tau)= B_kt^k$ where $B_k$ is the Bergman kernel of $L^k$, $\tau$ and $t$ are the coordinates on the second factor of $S^{2N_k+1}\times\R^+$ and $M\times\R^+$ respectively. 
	\end{prop}
	The same construction can be performed when the Sasakian manifold $M$ is the unitary bundle associated to the positive Hermitian bundle $(L,h)$ on a noncompact \K manifold $(X,\omega)$ with $\omega=-i\de\deb\log h$.
	Since the construction is similar we only highlight the necessary adjustments, see \cite{Rawnsley77CoherentStates} for details in the \K case.
	Firstly, we have to replace the vector space $H^0(L^k)$ with the Hilbert space of holomorphic sections of $L^k$ with finite norm
	\begin{equation}\label{EqHilbertSpace}
		\mathcal{H}_{k,h}=\left\{s\in \mathrm{Hol}(L) \ \vert\ \langle s,s\rangle_{k}<\infty\right\}.
	\end{equation}
	Now the construction makes sense only if we assume that the $\varepsilon$-function  
	\begin{equation}\label{EqReproducingKernel}
		\varepsilon_k(x)=\sum_{j=0}^N h^k\left( s_j(x),s_j(x)\right)
	\end{equation}
	of $\mathcal{H}_{k,h}$ with $N\leq\infty$ is a strictly positive function.
	Under this assumption the map $\varphi_k$ is well defined but the target space is generally $\ell^2(\C)\setminus\{0\}$ because the space of sections need not be finite dimensional. 
	In light of this, the map is now just a holomorphic immersion and the function $B_k$ is replaced by the $\varepsilon$-function $\varepsilon_k$.
	We sum up this discussion of the noncompact case in the following
	\begin{prop}\label{PropNoNCompactEmbedding}
		Let $M$ be the regular Sasakian manifold determined by the Hermitian bundle $(L,h)$ over a noncompact \K manifold $X$ and assume the space $\mathcal{H}_{k,h}$ is nontrivial. 
		Then there exists a holomorphic immersion $\varphi_k:M\times\R^+\lra S^{\infty}\times\R^+$ such that $\varphi_k^*(\tau)= \varepsilon_k t^k$ where $\varepsilon_k$ is the $\varepsilon$-function of $\mathcal{H}_{k,h}$, $\tau$ and $t$ are the coordinates on the second factor of $S^{\infty}\times\R^+$ and $M\times\R^+$ respectively. 
	\end{prop}
	
	\begin{remark}\rm\label{RmkPullbackCoordinate}
		Although $\varphi_k^*(h^{-1}_{FS})$ is not a hermitian metric on the line bundle $L^{-1}$ (it does not scale correctly under the $\C^*$-action), it defines a change of coordinate $(p,t)\mapsto (p,B_kt^k)$ (or $(p,t)\mapsto (p,\varepsilon_kt^k)$ in the noncompact case) on $M\times \R^+$ corresponding to the composition of the $\D_k$-homothetic transformation with a transverse \K deformation.
	\end{remark}

	\section{Sasakian immersions of homogeneous Sasakian manifolds into Sasakian space forms with $c>-3$}\label{SecImmersionPositive}
	Given the discussion of the previous section we are now ready to characterise embedding of homogeneous Sasakian manifolds into spheres.
	\begin{proof}[Proof of Theorem~\ref{TheoMainCpt}]
		Let $M$ be a compact homogeneous Sasakian manifold. 
		Then $M$ is the total space of the unit bundle associated to a holomorphic line bundle $\pi:L^{-1}\lra X$  over a generalised flag manifold $(X,\omega)$. 
		After possibly performing a $\D$-homothety, we can assume that $2\pi^*\omega=\di\eta$ is the transverse \K structure on $\D$ so that $c_1(L)=[\omega]$ where $L$ is the dual to $L^{-1}$.
		In particular, $L$ is an ample line bundle on a generalised flag manifold, hence very ample. 
		Thus, choosing a Hermitian metric $h$ on $L$ whose Ricci curvature form is $\omega$ and whose dual metric $h^{-1}$ on $L^{-1}$ induces the coordinate on $M\times \R^+$, we can apply Proposition~\ref{PropCompactEmbedding} with $k=1$.
		In particular, we get a a holomorphic embedding $\varphi=\varphi_1:M\times\R^+\lra S^{2N+1}\times\R^+$ such that $\varphi^*(\tau)= Bt$ where $B$ is the Bergman kernel of $L$, $\tau$ and $t$ are the coordinates on the second factor of $S^{2N+1}\times\R^+$ and $M\times\R^+$ respectively. 
		Moreover, since $(X,\omega)$ is a generalised flag manifold, the Bergman kernel $B$ is a positive constant, see \cite[Theorem~4.3]{ArezzoLoi04MomentMaps}.
		
		Now to get an isometry, i.e. to have $\varphi^*(\tau)=t$, it suffices to define the map $\varphi$ with respect to a rescaled orthogonal basis of $H^0(L)$.
		Namely, we consider a new orthogonal basis $(\widetilde s_0, \ldots, \widetilde s_N)$ of $H^0(L)$ by setting $\widetilde s_j=\dfrac{s_j}{\sqrt{B}}$ for $j=0,\ldots,N$. 
		It is immediate to compute the pullback of the Fubini-Study metric $h^{-1}_{FS}$ on $\OO(-1)=\C^{2N+1}$ as 
		\begin{align*}
			\varphi^*(h^{-1}_{FS})&=\sum_{j=0}^N\left\vert\alpha_v \widetilde f_j(x)\right\vert^2=\sum_{j=0}^N\left\vert v(\widetilde s_j(x))\right\vert^2=\dfrac{1}{B}\sum_{j=0}^N h^{-1}(v,v)h\left(s_j(x),s_j(x)\right)\\
			&=h^{-1}(v,v)\dfrac{\sum_{j=0}^N h\left(s_j(x),s_j(x)\right)}{B}=h^{-1}
		\end{align*}
		where we used the notation from Section~\ref{SecGeneralSasakianEmbedding}.
		
		Thus we have constructed a \K embedding $\varphi:M\times\R^+\lra S^{2N+1}\times\R^+$ of the \K cone of $M$ into the \K cone $\C^{N+1}\setminus\{0\}$ of the standard Sasakian sphere. 
		That is, a Sasakian embedding of $M$ into $S^{2N+1}$.
		
		Finally, let  $\varphi_1,\,\varphi_2$ be two \K embeddings of $M\times\R^+$ into $\C^{N+1}$ such that $\varphi_1^*(\tau)=\varphi_2^*(\tau)=t$ where $\tau$ is the Euclidean norm of $\C^{N+1}$. 
		By Calabi's rigidity there exists a rigid transformation $T$ of $\C^{N+1}$ such that $\varphi_1=T\circ\varphi_2$.
		Moreover, $T$ is unitary because $\varphi_1^*(\tau)=\varphi_2^*(\tau)=t$.
		We conclude that the embedding $\varphi$ is unique up to unitary transformations of $\C^{N+1}$.
	\end{proof}
	
	\begin{remark}\rm
		A compact homogeneous Sasakian manifold $M$ cannot admit an immersion at a point in any other Sasakian space form $M(n,c)$ for $c\leq-3$.
		This is because the generalised flag manifold $X$ over which $M$ fibres cannot be (locally) immersed in $F(N,b)$ when $b\leq 0$.
	\end{remark}
	
	We can now focus on the infinite dimensional case and classify homogeneous Sasakian manifolds which admit an immersion in $S^\infty$.
	
	\begin{proof}[Proof of~\Cref{TheoMainTop}]
		Let us first show that the fundamental group of a (not necessarily compact) homogeneous Sasakian manifold $M$ admitting an immersion in $S^{\infty}$ is cyclic. 
		Such a manifold $M$ fibres over a homogeneous \K manifold $X$ with fibre $\R$ or $S^1$.
		Moreover, the Sasakian immersion $M\lra S^{\infty}$ induces a \K immersion $X\lra\CP^\infty$.
		Hence $X$ is simply connected by \cite[Theorem~3]{Discala12Immersions} and the long exact sequence of the fibration yields the claim.
		
		Viceversa, let $M$ be a homogeneous Sasakian manifold with cyclic fundamental group.
		If $M$ is compact, then the claim follows directly from Theorem~\ref{TheoMainCpt}.
		So we can restrict to the case where $M$ is noncompact.
		By Proposition~\ref{PropSasakiNonCompactfibreR}, $M$ fibres over a noncompact homogeneous \K manifold $(X,\omega)$. 
		By a classical result of Dorfmeister and Nakajima \cite{Dorfmeister88ClassificationHomogeneous}, such a manifold $X$ is a complex (but not necessarily K\"ahler) product $\C^d\times T\times F\times \Omega$ where $T$ is a complex torus, $F$ is a generalised flag manifold, and $\Omega$ a homogeneous bounded domain.
		From the long exact sequence of the fibration we deduce that the factor $T$ cannot appear because $\pi_1(X)=\pi_1(T)$ must be the quotient of the cyclic group $\pi_1(M)$.
		Therefore $X$ is simply connected.

		Let us first treat the case where the fibre of the bundle $M\lra X$ is $S^1$ and deal later with the case in which the fibre is $\R$.
		Although $M$ is noncompact, the Sasakian structure on $M$ is regular and the fibre is $S^1$ so that, analogously to the compact case, $M$ is the unit bundle of a line bundle $\pi:L^{-1}\lra X$. Moreover, after possibly performing a $\D$-homothety, we can assume that $2\pi^*\omega=\di\eta$ is the transverse \K structure on $\D$ so that $c_1(L)=[\omega]$ where $L$ is the dual to $L^{-1}$.
		Choose a hermitian metric $h$ on $L$ whose Ricci curvature form is $\omega$. 
		Now there exists a suitable real $k>0$ such that the space $\mathcal{H}_{k,h}$ defined in Section~\ref{SecGeneralSasakianEmbedding} is nontrivial, see \cite{rosenberg85Harmonically}.
		Moreover, $\mathcal{H}_{k,h}\neq\{0\}$ if and only if its $\varepsilon$-function $\varepsilon_k$ is in fact a positive constant, cf. \cite[Lemma~2.2]{Loi15HomogeneousImmersions}.
		Thus, for such a $k$, the hypotheses of Proposition~\ref{PropNoNCompactEmbedding} are satisfied.
		We then get a holomorphic immersion $\varphi_k:M\times\R^+\lra S^{\infty}\times\R^+$ such that
		$\varphi_k^*(\tau)= \varepsilon_k t^k$ where $\varepsilon_k$ is constant, $\tau$ and $t$ are the coordinates on the second factor of $S^{\infty}\times\R^+$ and $M\times\R^+$ respectively. 
		Similarly to the proof of Theorem~\ref{TheoMainCpt}, by rescaling the orthogonal basis used in the construction of the immersion $\varphi_k$, we can modify the immersion so that $\varphi_k^*(\tau)=t^k$.
		Recall now that replacing $t$ by $t^k$ is equivalent to taking a $\D_k$-homothety (cf. Remark~\ref{RmkPullbackCoordinate} and Definition\ref{DefinTransverseHomothety}) so that we have proven the statement when the fibre is $S^1$.

		We are left now with the case where the fibre of the fibration $M\lra X$ is $\R$.
		By Proposition~\ref{PropSasakiNonCompactfibreR}, $M$ fibres over a contractible homogeneous \K manifold $(X,\omega)$, that is, $M=X\times\R$.
		Moreover, $X$ is a homogeneous \K manifold with cyclic, hence trivial, fundamental group.
		Notice that $\Z\subset\R$ acts on $X\times\R$ by Sasakian isometries via the flow of the Reeb vector field. 
		The quotient is the Sasakian manifold $N=X\times S^1$ and the universal covering map $\widetilde \pi:M\lra N$ is a Sasakian immersion (it is simply the exponential map on the second factor).
		Now $N$ is a Sasakian manifold fibring over a simply connected homogeneous \K manifold $X$.
		By the previous case, there exists a Sasakian immersion $\varphi:N\lra S^\infty$ of a suitable $\D$-homothety of $N$ into $S^\infty$. 
		The pullback of this structure to $M$ under $\widetilde\pi$ is a $\D$-homothetic transformation of the homogeneous structure we began with.
		Finally, the composition $\widetilde \pi\circ\varphi:M\lra S^\infty$ is a Sasakian immersion of $M$ (endowed with such a $\D$-homothetic structure) into $S^\infty$. 
	\end{proof}

	\begin{remark}\label{necessaryDhomothetic}\rm
		The assumption on the $\D$-homothetic transformation is necessary because there exist simply connected homogeneous Sasakian manifolds which do not admit an immersion into $S^\infty$ unless a transverse homothety is performed.
		Such an example is given by a homogeneous Sasakian structure on $M=\Omega\times\R$ where the transverse \K structure is given by a suitable multiple $c\omega_B$ of the Bergman metric on a bounded symmetric domain $\Omega$. 
		In fact, a Sasakian immersion $M\lra S^\infty$ covers a \K immersion $\Omega\lra\CP^\infty$ and this exists if and only if $c\gamma$ belongs to the Wallach set of $\Omega$ (see \cite[Theorem~2]{Loi11InfiniteProjective}), where $\gamma$ denotes the genus of $\Omega$. Therefore, any bounded symmetric domain $\Omega$ (other than $\CH^n$) and small enough $c$ provide such an instance.
	\end{remark}
	
	\begin{remark}\label{infinitenoncompact}\rm
		A homogeneous Sasakian manifold $M$ which admits a Sasakian immersion in $S^{2N+1}$ with $N<\infty$ is necessarily compact. 
		In fact it fibres over a homogeneous \K manifold $X$ which admits a \K immersion in $\CP^N$. 
		By \cite{Takeuchi78Homogeneous} $X$ is forced to be compact and the fibre of $M\lra X$ is forced to be $S^1$.
		This implies that the immersion in Theorem~\ref{TheoMainTop} is full, i.e. its image is not contained in a proper totally geodesic submanifold, if and only if the manifold $M$ is noncompact.
		This is no longer true if we drop the homogeneity assumption, even if the metric is complete and $N<\infty$. Such an immersion is given in the following example.
	\end{remark}
	
	\begin{ex}\label{ExCompleteFiniteInjective}\rm
		Let $T\subset\CP^n$ be a complex $2$-dimensional abelian variety and choose a copy of $\C$ dense in $T$. Restrict the Fubini-Study metric of $\CP^n$ to $T$ and to $\C$ denoting them with $g_T$ and $g$ respectively. 
		Thus $(\C,g)$ is a noncompact \K manifold  isometrically immersed in $\CP^n$.
		Now the Hopf bundle $S^{2n+1}\lra \CP^n$ restricts to a circle bundle $M\lra\C$.
		The manifold $M$ is a noncompact regular Sasakian manifold which admits an injective Sasakian immersion into $S^{2n+1}$, hence $S^\infty$.
		Moreover, the Sasakian metric on $M$ is complete because the \K metric $g$ on $\C$ is.
		In fact, by compactness, the \K metric $g_T$ satisfies $ag_0<g_T<bg_0$ where $g_0$ is the flat metric and $a,b$ are suitable constants.
		Therefore, the same inequalities hold on $\C$, i.e. $ag_0<g<bg_0$, forcing $g$ to be complete.
	\end{ex}

	\bibliographystyle{amsplain}
	\bibliography{biblio}

\providecommand{\bysame}{\leavevmode\hbox to3em{\hrulefill}\thinspace}
\providecommand{\MR}{\relax\ifhmode\unskip\space\fi MR }
\providecommand{\MRhref}[2]{%
  \href{http://www.ams.org/mathscinet-getitem?mr=#1}{#2}
}
\providecommand{\href}[2]{#2}
\begin{thebibliography}{10}

\bibitem{ArezzoLoi04MomentMaps}
Claudio Arezzo and Andrea Loi, \emph{Moment maps, scalar curvature and
  quantization of {K}\"{a}hler manifolds}, Comm. Math. Phys. \textbf{246}
  (2004), no.~3, 543--559. \MR{2053943}

\bibitem{Bande20EtaEinsteinNonCompact}
Gianluca Bande, Beniamino Cappelletti-Montano, and Andrea Loi,
  \emph{{$\eta$}-{E}instein {S}asakian immersions in non-compact {S}asakian
  space forms}, Ann. Mat. Pura Appl. (4) \textbf{199} (2020), no.~6,
  2117--2124. \MR{4165673}

\bibitem{Boyer08Book}
Charles~P. Boyer and Krzysztof Galicki, \emph{Sasakian geometry}, Oxford
  Mathematical Monographs, Oxford University Press, Oxford, 2008. \MR{2382957}

\bibitem{Calabi53Isometric}
Eugenio Calabi, \emph{Isometric imbedding of complex manifolds}, Ann. of Math.
  (2) \textbf{58} (1953), 1--23. \MR{57000}

\bibitem{Cappelletti19EinsteinSpheres}
Beniamino Cappelletti-Montano and Andrea Loi, \emph{Einstein and {$\eta
  $}-{E}instein {S}asakian submanifolds in spheres}, Ann. Mat. Pura Appl. (4)
  \textbf{198} (2019), no.~6, 2195--2205. \MR{4031847}

\bibitem{Discala12Immersions}
Antonio~Jose Di~Scala, Hideyuki Ishi, and Andrea Loi, \emph{K\"{a}hler
  immersions of homogeneous {K}\"{a}hler manifolds into complex space forms},
  Asian J. Math. \textbf{16} (2012), no.~3, 479--487. \MR{2989231}

\bibitem{Dorfmeister88ClassificationHomogeneous}
Josef Dorfmeister and Kazufumi Nakajima, \emph{The fundamental conjecture for
  homogeneous {K}\"{a}hler manifolds}, Acta Math. \textbf{161} (1988), no.~1-2,
  23--70. \MR{962095}

\bibitem{Green78NonRigidity}
Mark~L. Green, \emph{Metric rigidity of holomorphic maps to {K}\"{a}hler
  manifolds}, J. Differential Geometry \textbf{13} (1978), no.~2, 279--286.
  \MR{540947}

\bibitem{Harada72SpaceFormInSpaceForm}
Minoru Harada, \emph{Sasakian space forms immersed in {S}asakian space forms},
  Bull. Tokyo Gakugei Univ. (4) \textbf{24} (1972), 7--11. \MR{322745}

\bibitem{Harada72SubmanifoldsI}
\bysame, \emph{On {S}asakian submanifolds}, Tohoku Math. J. (2) \textbf{25}
  (1973), 103--109. \MR{324590}

\bibitem{Harada72SubmanifoldsII}
\bysame, \emph{On {S}asakian submanifolds. {II}}, Bull. Tokyo Gakugei Univ. (4)
  \textbf{25} (1973), 19--23. \MR{324591}

\bibitem{Harada75SomeSasakianSubmanifolds}
\bysame, \emph{On some {S}asakian submanifolds immersed in {S}asakian space
  forms}, Bull. Tokyo Gakugei Univ. (4) \textbf{27} (1975), 35--39. \MR{388277}

\bibitem{Kenmotsu69InvariantSubmanifolds}
Katsuei Kenmotsu, \emph{Invariant submanifolds in a {S}asakian manifold},
  Tohoku Math. J. (2) \textbf{21} (1969), 495--500. \MR{259806}

\bibitem{Kon73Kodai}
Masahiro Kon, \emph{Invariant submanifolds of normal contact metric manifolds},
  K\={o}dai Math. Sem. Rep. \textbf{25} (1973), 330--336. \MR{324581}

\bibitem{Kon76InvariantSubmanifolds}
\bysame, \emph{Invariant submanifolds in {S}asakian manifolds}, Math. Ann.
  \textbf{219} (1976), no.~3, 277--290. \MR{425844}

\bibitem{Loi15HomogeneousImmersions}
Andrea Loi and Roberto Mossa, \emph{Some remarks on homogeneous {K}\"{a}hler
  manifolds}, Geom. Dedicata \textbf{179} (2015), 377--383. \MR{3424675}

\bibitem{Loi11InfiniteProjective}
Andrea Loi and Michela Zedda, \emph{K\"{a}hler-{E}instein submanifolds of the
  infinite dimensional projective space}, Math. Ann. \textbf{350} (2011),
  no.~1, 145--154. \MR{2785765}

\bibitem{Loi18Book}
\bysame, \emph{K\"{a}hler immersions of {K}\"{a}hler manifolds into complex
  space forms}, Lecture Notes of the Unione Matematica Italiana, vol.~23,
  Springer, Cham; Unione Matematica Italiana, [Bologna], 2018. \MR{3838438}

\bibitem{Okumura68Immersion}
Masafumi Okumura, \emph{On contact metric immersion}, K\={o}dai Math. Sem. Rep.
  \textbf{20} (1968), 389--409. \MR{264564}

\bibitem{Ornea07CREmbedding}
Liviu Ornea and Misha Verbitsky, \emph{Embeddings of compact {S}asakian
  manifolds}, Math. Res. Lett. \textbf{14} (2007), no.~4, 703--710.
  \MR{2335996}

\bibitem{Placini21SasakiSolitons}
G.~Placini, \emph{Sasakian immersions of {S}asaki-{R}icci solitons into
  {S}asakian space forms}, J. Geom. Phys. \textbf{166} (2021), Paper No.
  104265, 7. \MR{4249127}

\bibitem{Rawnsley77CoherentStates}
J.~H. Rawnsley, \emph{Coherent states and {K}\"{a}hler manifolds}, Quart. J.
  Math. Oxford Ser. (2) \textbf{28} (1977), no.~112, 403--415. \MR{466649}

\bibitem{Reinhart59Foliated}
Bruce~L. Reinhart, \emph{Foliated manifolds with bundle-like metrics}, Ann. of
  Math. (2) \textbf{69} (1959), 119--132. \MR{107279}

\bibitem{rosenberg85Harmonically}
Jonathan Rosenberg and Mich\`ele Vergne, \emph{Harmonically induced
  representations of solvable {L}ie groups}, J. Funct. Anal. \textbf{62}
  (1985), no.~1, 8--37. \MR{790768}

\bibitem{Takeuchi78Homogeneous}
Masaru Takeuchi, \emph{Homogeneous {K}\"{a}hler submanifolds in complex
  projective spaces}, Japan. J. Math. (N.S.) \textbf{4} (1978), no.~1,
  171--219. \MR{528871}

\bibitem{Tanno69SpaceForm}
Sh\^{u}kichi Tanno, \emph{Sasakian manifolds with constant {$\phi$}-holomorphic
  sectional curvature}, Tohoku Math. J. (2) \textbf{21} (1969), 501--507.
  \MR{251667}

\end{thebibliography}
\end{document}